# SIEGEL MODULAR FORMS AND THETA SERIES ATTACHED TO QUATERNION ALGEBRAS II

SIEGFRIED BOECHERER AND RAINER SCHULZE-PILLOT

ABSTRACT. We continue our study of Yoshida's lifting, which associates to a pair of automorphic forms on the adelic multiplicative group of a quaternion algebra a Siegel modular form of degree 2. We consider here the case that the automorphic forms on the quaternion algebra correspond to modular forms of arbitrary even weights and square free levels; in particular we obtain a construction of Siegel modular forms of weight 3 attached to a pair of elliptic modular forms of weights 2 and 4.

## INTRODUCTION

We resume in this article our study of Yoshida's lifting from [5]. This lifting associates a Siegel modular form $Y^{(2)}(\varphi_1, \varphi_2)$ of degree 2 to a pair $\varphi_1, \varphi_2$ of automorphic forms on the adelic multiplicative group of a definite quaternion algebra $D$ over $\mathbf{Q}$. This pair corresponds under Eichler's correspondence to a pair $f, g$ of elliptic modular forms of weights $k_1 = 2, k_2$ and the same square free level $N$; it is a theta lifting from the orthogonal group of $D$ equipped with the norm form to the group $Sp_2$. In the case that also $k_2 = 2$ holds we proved that this lifting is nonzero if both $f, g$ are cusp forms with the same eigenvalues under the Atkin-Lehner involutions; if $\varphi_2$ is constant (and hence $g$ is an Eisenstein series) the vanishing of $Y^{(2)}(\varphi_1, 1)$ depends on the central critical value of the $L$-function of $f$. We extended the construction to arbitrary pairs of weights $k_1, k_2$ in [6], where we also gave the analogous construction of a lifting $Y^{(n)}(\varphi_1, \varphi_2)$ to Siegel modular forms of higher degree. These higher degree liftings are in general vector valued Siegel modular forms; $Y^{(n)}(\varphi_1, \varphi_2)$ is mapped to $Y^{(n-1)}(\varphi_1, \varphi_2)$ by Siegel's $\phi$-operator.

The goal of this article is to extend our nonvanishing results from [5] to this more general situation. The method of proof is essentially the same: To prove the nonvanishing of $Y^{(2)}(\varphi_1, \varphi_2)$ we show that $Y^{(3)}(\varphi_1, \varphi_2)$ is not cuspidal. This is done by studying the analytic behaviour of its standard $L$-function in $s = 1$ in two ways: Once

Both authors supported by MSRI, Berkeley (NSF-grant DMS-9022140). R. Schulze-Pillot also supported by the Deutsche Forschungsgemeinschaft and by the Max-Planck-Institut für Mathematik, Bonn.

by giving an integral representation (assuming $Y^{(3)}(\varphi_1, \varphi_2)$ to be cuspidal) and once by computing the $L$-function in terms of the Satake parameters of $f, g$ using the construction of $Y^{(3)}(\varphi_1, \varphi_2)$ as a theta lifting. The integral representation which we use is the same as in [5] using the doubling (or pullback) method and an Eisenstein series of weight 2, modified by some differential operator. In the generality needed for our purpose these differential operators have been introduced and studied by Ibukiyama [17]. On the theta lifting side we arrive after a careful computation of the factors occurring at the bad places at the study of the value at $s = 1$ of the normalized $GL_2 \times GL_2 - L-$ function (with functional equation under $s \mapsto 1 - s$) associated to $f, g$. A result of Ogg [23] and Shahidi [24] gives the nonvanishing of this value if $f, g$ are cuspidal and nonproportional, which contradicts the results obtained from the integral representation (and hence our assumption that $Y^{(3)}(\varphi_1, \varphi_2)$ is cuspidal). The case where $\varphi_2$ is constant is more delicate and leads to the $L$-function condition mentioned.

To summarize our results, let $f, g$ be primitive cusp forms of square free levels $N_f, N_g$ (with $N = \mathrm{lcm}(N_f, N_g), \gcd(N_f, N_g) \neq 1$), weights $k_1 \geq k_2$ and Haupttyp for $\Gamma_0(N_{f,g})$ (eigenforms of all Hecke operators) and assume that for all $p$ dividing $\gcd(N_f, N_g)$ the eigenvalues of $f, g$ under the Atkin-Lehner-involution $w_p$ are the same. Then we can associate to $f, g$ a nonzero vector valued Siegel cusp form of degree 2, whose Satake parameters are computed from those of $f, g$. If $k_2 = 2$ this Siegel modular form is scalar valued of weight $1 + \frac{k_1}{2}$; in particular for $k_1 = 4, k_2 = 2$ we obtain cusp forms of weight 3. These are of interest in geometry since they correspond to holomorphic differential forms of top degree on the quotient of the Siegel upper half space of degree 2 by the group $\Gamma_0(N)$. It should be noted that we obtain in fact different constructions for each subset of odd cardinality of the set of prime divisors of $\gcd(N_f, N_g)$ and that the resulting Siegel modular forms are pairwise orthogonal. Most of this article was written while both authors were guests of the MSRI in Berkeley during its special year on automorphic forms. We wish to thank the MSRI for its hospitality and financial support. R. Schulze-Pillot was also supported by Deutsche Forschungsgemeinschaft during a visit of one month at MSRI and was a guest of the Max-Planck-Institut für Mathematik in Bonn in the final stage of the preparation of this manuscript.

## 1. Yoshida's lifting

Yoshida's lifting [30, 31] associates a Siegel modular form of degree 2 to a pair $f, g$ of elliptic modular forms of even weight, where one of $f, g$ has weight 2. Its generalization to arbitrary even weights considered in [6] yields vector valued Siegel modular forms. As explained in the introduction we want to generalize our results on the nonvanishing of Yoshida's lifting to this more general situation. For this, we need to recall a few notations and results from [5, 6].

Let $D$ be a definite quaternion algebra over $\mathbf{Q}$ and $R$ an Eichler order of square free level $N$ in $D$, i. e. , the completion $R_p$ is a maximal order in $D_p$ for the primes $p$ ramified in $D$ and is of index $p$ in a maximal order (i. e. $R_p$ is conjugate to the set of matrices $\begin{pmatrix} a & b \\ c & d \end{pmatrix}$ in $M_2(\mathbf{Z}_p)$ with $p \mid c$) for the remaining primes $p$ dividing $N$. We will usually decompose $N$ as $N = N_1 N_2$ where $N_1$ is the product of the primes that are ramified in $D$. On $D$ we have the involution $x \mapsto \bar{x}$, the (reduced) trace $\operatorname{tr}(x) = x + \bar{x}$ and the (reduced) norm $n(x) = x\bar{x}$.

For $\nu \in \mathbf{N}$ let $U_\nu^{(0)}$ be the space of homogeneous harmonic polynomials of degree $\nu$ on $\mathbf{R}^3$ and view $P \in U_\nu^{(0)}$ as a polynomial on $D_\infty^{(0)} = \{x \in D_\infty | \operatorname{tr}(x) = 0\}$ by putting $P(\sum_{i=1}^3 x_i e_i) = P(x_1, x_2, x_3)$ for an orthonormal basis $\{e_i\}$ of $D_\infty^{(0)}$ with respect to the norm form $n$.

The representations $\tau_\nu$ of $D_\infty^\times/\mathbf{R}^\times$ of highest weight $(\nu)$ on $U_\nu^{(0)}$ given by

$$(\tau_\nu(y))(P)(x) = P(y^{-1}xy)$$

for $\nu \in \mathbf{N}$ give all the isomorphism classes of irreducible rational representations of $D_\infty^\times/\mathbf{R}^\times$. By $\langle\!\langle \ , \ \rangle\!\rangle$ we denote the suitably normalized invariant scalar product in the representation space $U_\nu^{(0)}$.

The group of proper similitudes of the quadratic form $q(x) = n(x)$ (with associated symmetric bilinear form $B(x,y) = \operatorname{tr}(x\bar{y})$) on $D$ is isomorphic to $(D^\times \times D^\times)/Z(D^\times)$ (as algebraic group) via

$$(x_1, x_2) \mapsto \sigma_{x_1,x_2} \text{ with } \sigma_{x_1,x_2}(y) = x_1 y x_2^{-1},$$

with the special orthogonal group being the image of

$$\{(x_1, x_2) \in D^\times \times D^\times \mid n(x_1) = n(x_2)\}.$$

We denote by $H$ the orthogonal group of $(D, n)$, by $H^+$ the special orthogonal group and by $K$ (resp. $K^+$) the group of units (resp. proper units) of the lattice $R$ in $D$. Let $\nu_1 \geq \nu_2$ be given and let $\tau_1 = \tau_{\nu_1}, \tau_2 = \tau_{\nu_2}$ for $i = 1, 2$. It is then well known that the $H^+(\mathbf{R})$-space $U_{\nu_1}^{(0)} \otimes U_{\nu_2}^{(0)}$ is isomorphic to the $H^+(\mathbf{R})$-space $U_{\nu_1,\nu_2}$ of $\mathbf{C}[X_1, X_2]$-valued harmonic forms on $D_\infty^2$ transforming according to the representation of $GL_2(\mathbf{R})$ of highest weight $(\nu_1 + \nu_2, \nu_1 - \nu_2)$; an intertwining map $\Psi$ has been given in [6, Section 3] (for the theory of harmonic forms see [10, 19]). It is also well known [19] that the representation $\lambda_{\nu_1,\nu_2}$ of $H^+(\mathbf{R})$ on $U_{\nu_1,\nu_2}$ is irreducible of highest weight $(\nu_1 + \nu_2, \nu_1 - \nu_2)$. If $\nu_1 \neq \nu_2$ it can be extended in a unique way to an irreducible representation of $H(\mathbf{R})$ on the space $U_{\nu_1,\nu_2,s} := (U_{\nu_1}^{(0)} \otimes U_{\nu_2}^{(0)}) \oplus (U_{\nu_2}^{(0)} \otimes U_{\nu_1}^{(0)})$, denoted by $(\tau_1 \otimes \tau_2)_s$, whereas for $\nu_1 = \nu_2$ there are two possible extensions to representations $(\tau_1 \otimes \tau_2)_\pm$ on $U_{\nu_1,\nu_2}$; we denote this space with the representation $(\tau_1 \otimes \tau_2)_\pm$ on it by $U_{\nu_1,\nu_2,\pm}$.

For an irreducible rational representation $(V_\tau, \tau)$ (with $\tau = \tau_\nu$ as above) of $D_\infty^\times/\mathbf{R}^\times$ we denote by $\mathcal{A}(D_\mathbf{A}^\times, R_\mathbf{A}^\times, \tau)$ the space of functions $\varphi : D_\mathbf{A}^\times \to V_\tau$ satisfying $\varphi(\gamma x u) =$ $\tau(u_\infty^{-1})\varphi(x)$ for $\gamma \in D_{\mathbf{Q}}^\times$ and $u = u_\infty u_f \in R_{\mathbf{A}}^\times$, where $R_{\mathbf{A}}^\times = D_\infty^\times \times \prod_p R_p^\times$ is the adelic group of units of $R$.

It has been discovered by Eichler that these functions are in correspondence with the elliptic modular forms of weight $2 + 2\nu$ and level $N$. This correspondence can be described as follows:

Let $D_{\mathbf{A}}^\times = \cup_{i=1}^r D^\times y_i R_{\mathbf{A}}^\times$ be a double coset decomposition with $y_{i,\infty} = 1$ and $n(y_i) = 1$, put $I_{ij} = y_i R y_j^{-1}$, $R_i = I_{ii}$ and let $e_i$ be the number of units of the order $R_i$.

Recall from §5 of [5] that for each $p \mid N$ we have an element $\pi_p \in D_p^\times$ of norm $p$ normalizing $R_p$ and that these $\pi_p$ (together with $R_{\mathbf{A},f}^\times$) generate the normalizer of $R_{\mathbf{A},f}^\times \mathbf{Q}_{\mathbf{A}}^\times / \mathbf{Q}_{\mathbf{A}}^\times$ in $D_{\mathbf{A},f}^\times \mathbf{Q}_{\mathbf{A}}^\times / \mathbf{Q}_{\mathbf{A}}^\times$ and satisfy $\pi_p^2 \in \mathbf{Q}_p R_p^\times$. Right translation by $\pi_p$ gives an involution $\widetilde{w_p}$ of $\mathcal{A}(D_{\mathbf{A}}^\times, R_{\mathbf{A}}^\times, \tau)$, which space then splits into common eigenspaces of all these (pairwise commuting) involutions.

On the space $\mathcal{A}(D_{\mathbf{A}}^\times, R_{\mathbf{A}}^\times, \tau)$ we have for $p \nmid N$ Hecke operators $\tilde{T}(p)$ defined by

$$\tilde{T}(p)\varphi(x) = \int_{D_p^\times} \varphi(xy^{-1})\chi_p(y)\, dy$$

where $\chi_p$ is the characteristic function of $\{y \in R_p | n(y) \in p\mathbf{Z}_p^\times\}$. They commute with the involutions $\tilde{w}_p$ and are given explicitly by $\tilde{T}(p)\varphi(y_i) = \sum_{j=1}^r B_{ij}(p)\varphi(y_j)$, where the Brandt matrix entry $B_{ij}(p)$ is given as

$$B_{ij}(p) = B_{ij}^{(\nu)}(p) = \frac{1}{e_j} \sum_{\substack{x \in y_j R y_i^{-1} \\ n(x) = p}} p^\nu \tau(x) ,$$

hence is itself an endomorphism of the representation space $U_\nu^{(0)}$ of $\tau$. On the space $\mathcal{A}(D_{\mathbf{A}}^\times, R_{\mathbf{A}}^\times, \tau)$ we have the natural inner product $\langle\ ,\ \rangle$ defined by integration, it is explicitly given by

$$\langle \varphi, \psi \rangle = \sum_{i=1}^h \frac{\langle\langle \varphi(y_i), \psi(y_i) \rangle\rangle}{e_i}.$$

By abuse of language we call (in the case $\nu = 0$) forms cuspidal, if they are orthogonal to the constant functions with respect to this inner product.

From [8, 14, 25, 18] we know then that the essential part $\mathcal{A}_{\text{ess}}(D_{\mathbf{A}}^\times, R_{\mathbf{A}}^\times, \tau)$ consisting of functions $\varphi$ that are orthogonal to all $\psi \in \mathcal{A}(D_{\mathbf{A}}^\times, (R_{\mathbf{A}}')^\times, \tau)$ for orders $R'$ strictly containing $R$ is invariant under the $\tilde{T}(p)$ for $p \nmid N$ and the $\tilde{w}_p$ for $p \nmid N$ and hence has a basis of common eigenfunctions of all the $\tilde{T}(p)$ for $p \nmid N$ and all the involutions $\tilde{w}_p$ for $p \mid N$. Moreover in $\mathcal{A}_{\text{ess}}(D_{\mathbf{A}}^\times, R_{\mathbf{A}}^\times, \tau)$ strong multiplicity one holds, i.e., each system of eigenvalues of the $\tilde{T}(p)$ for $p \nmid N$ occurs at most once, and the eigenfunctions are in one to one correspondence with the newforms in the space $S^{2+2\nu}(N)$ of elliptic cusp forms of weight $2 + 2\nu$ for the group $\Gamma_0(N)$ that are eigenfunctions of all Hecke operators (if $\tau$ is the trivial representation and $R$ is a maximal order one has to restrict here to functions orthogonal to the constant function 1 on the quaternion side in order to obtain cusp forms on the modular forms side). This correspondence (Eichler's correspondence) preserves Hecke eigenvalues for $p \nmid N$, and if $\varphi$ corresponds to $f \in S^{2+2\nu}(N)$ then the eigenvalue of $f$ under the Atkin-Lehner involution $w_p$ is equal to that of $\varphi$ under $\tilde{w}_p$ if $D$ splits at $p$ and equal to minus that of $\varphi$ under $\tilde{w}_p$ if $D_p$ is a skew field. An extension of this correspondence to forms $\varphi$ as above that are not essential but eigenfunctions of all the involutions $\tilde{w}_p$ has been given in [13, 7].

Lemma 3.1 of [6] now gives the following setup for the Yoshida lifting:

Let $\mathcal{A}(H_{\mathbf{A}}^+, K_{\mathbf{A}}^+, \tau_1 \otimes \tau_2)$ be defined in the same way as $A(D_{\mathbf{A}}^\times, R_{\mathbf{A}}^\times, \tau_i)$. If $\varphi_i \in \mathcal{A}(D_{\mathbf{A}}^\times, R_{\mathbf{A}}^\times, \tau_i)$ for $i = 1, 2$ (with $\tau_i = \tau_{\nu_i}$ and $\nu_1 \geq \nu_2$) are eigenfunctions of all the involutions $\tilde{w}_p$ for $p \mid N$ and of all the Hecke operators $\tilde{T}_p$ for $p \nmid N$, we denote by $\widetilde{\varphi_1 \otimes \varphi_2}$ the function on the adelic special orthogonal group $H_{\mathbf{A}}^+$ of $(D, n)$ given by

$$\widetilde{\varphi_1 \otimes \varphi_2}(\sigma_{y_1,y_2}) = \int_{K_f^+} \varphi_1(y_1 k_1) \otimes \varphi_2(y_2 k_2) dk,$$

where $k$ in the finite part $K_f^+$ of the adelic group of proper units of the lattice $R$ in $D$ is represented as $k = \sigma_{k_1,k_2}$. Then $\widetilde{\varphi_1 \otimes \varphi_2} \neq 0$ if and only if $\varphi_1$ and $\varphi_2$ have the same eigenvalues under all the involutions $\tilde{w}_p$ for the $p \mid N$, and one has $\widetilde{\varphi_1 \otimes \varphi_2} \in \mathcal{A}(H_{\mathbf{A}}^+, K_{\mathbf{A}}^+, \tau_1 \otimes \tau_2)$. Let now $\varphi_1$ and $\varphi_2$ be Hecke eigenforms having the same eigenvalues under all the involutions $\tilde{w}_p$ for the $p \mid N$. If $(\tau_1 \otimes \tau_2)_*$ for $*$ denoting $+, -$ or $s$ denotes the possible extensions of $\tau_1 \otimes \tau_2$ to an irreducible representation of $H_{\mathbf{R}}$ let $(\varphi_1 \otimes \varphi_2)_*$ be the unique function in $\mathcal{A}(H_{\mathbf{A}}, K_{\mathbf{A}}, (\tau_1 \otimes \tau_2)_*)$ whose value at $h \in H_{\mathbf{A}}^+$ is

$$\frac{1}{2}(\widetilde{\varphi_1 \otimes \varphi_2}(h) + (\tau_1 \otimes \tau_2)_*(\iota)\widetilde{\varphi_1 \otimes \varphi_2}(\iota h \iota)),$$

where $\iota \in H_{\mathbf{R}}$ is the involution in $D$.

Then for $\nu_1 \neq \nu_2$ one has $(\varphi_1 \otimes \varphi_2)_s \neq 0$ if and only if $\widetilde{\varphi_1 \otimes \varphi_2} \neq 0$. For $\nu_1 = \nu_2$ one has $(\varphi_1 \otimes \varphi_2)_+ \neq 0$. Moreover in this case one has $(\varphi_1 \otimes \varphi_2)_- = 0$ if and only if $\varphi_1 = 0$ or $\varphi_2 = a\varphi_1$ for some $a \in \mathbf{C}$.

The functions $(\varphi_1 \otimes \varphi_2)_*$ on the adelic orthogonal group $H_{\mathbf{A}}$ thus constructed can be lifted to Siegel modular forms of degree $n$, we call this lifting Yoshida's lifting, denoted by $Y^{(n)}(\varphi_1, \varphi_2, *)$. To describe it explicitly we recall from [19, 29] that for any irreducible representation $(\lambda, U_\lambda)$ of $H(\mathbf{R})$ the space $\mathcal{H}_n(\lambda)$ of pluriharmonic polynomials $P : M_{m,n}(\mathbf{C}) \to U_\lambda$ such that $P(h^{-1}x) = \lambda(h^t)P(x)$ for all $h \in O_m$ is zero or (under the right action of $GL_n(\mathbf{C})$ on the variable) isomorphic to an irreducible representation $(\rho_n(\lambda), W_{\rho_n(\lambda)})$ of $GL_n(\mathbf{C})$. In the latter case the space $\mathcal{H}_q(\rho_n(\lambda))$ consisting of all $q$-pluriharmonic polynomials $P : M_{m,n}(\mathbf{C}) \to W_{\rho_n(\lambda)}$ such that $P(xg) = (\rho_n(\lambda)(g^t))P(x)$ for all $g \in GL_n(\mathbf{C})$ is isomorphic to $(U_\lambda, \lambda)$ as a representation space of $H_{\mathbf{R}}$.

We denote by $\mathcal{P}_{n,*}$ the (essentially unique) isomorphism from $U_{\nu_1,\nu_2,*}$ to $\mathcal{H}_q(\rho_n((\tau_1 \otimes \tau_2)_*))$.

Then, again for $*$ denoting $+$, $-$ or $s$, we have the $n$-th Yoshida-lifting (whenever the representation $\rho_n((\tau_1 \otimes \tau_2)_*)$ is defined):

$$\begin{aligned} Y^{(n)}(\varphi_1, \varphi_2, *) &:= \Theta^{(n)}(\varphi_1 \otimes \varphi_2)_*(Z) \\ &= \int_{H_{\mathbf{Q}}^+ \backslash H_{\mathbf{A}}^+} \sum_{\mathbf{x} \in (hL)^n} \mathcal{P}_{n,*}((\widetilde{\varphi_1 \otimes \varphi_2})(h))(h_\infty^{-1} x_1, \ldots, h_\infty^{-1} x_n) \times \\ &\quad \times \exp(2\pi i tr(q(\mathbf{x})Z)) dh \end{aligned}$$

It is explicitly given as

$$\begin{aligned} Y^{(n)}(\varphi_1, \varphi_2, *)(Z) &= \sum_{i,j=1}^r \frac{1}{e_i e_j} \sum_{\mathbf{x} \in (y_i R y_j^{-1})^n} \mathcal{P}_{n,*}(\Psi(\varphi_1(y_i) \otimes \varphi_2(y_j)))(x_1, \ldots, x_n) \times \\ &\quad \times \exp(2\pi i tr(q(\mathbf{x})Z)) \end{aligned}$$

In particular, if $\nu_1 = \nu_2 = \nu$ then $Y^{(1)}(\varphi_1, \varphi_2, +)$ is defined and is an elliptic modular form of weight $2 + 2\nu$. If $\varphi_1, \varphi_2$ are eigenforms of the $\tilde{T}_p$ for $p \nmid N$ with different eigenvalues then as in [31] one sees that $Y^{(1)}(\varphi_1, \varphi_2, +) = 0$. If $\varphi_1 = \varphi_2 = \varphi \neq 0$ is an eigenform of all the $\widetilde{T}_p$ for $p \nmid N$ and of all the $\widetilde{w}_p$ for the $p \mid N$ then $Y^{(1)}(\varphi, \varphi, +)$ is nonzero, and linear continuation of the map sending an eigenform $\varphi$ to $Y^{(1)}(\varphi, \varphi)$ realizes Eichler's correspondence described above. For details in the case that $\varphi$ is not essential see [7].

To investigate the nonvanishing of the lifting in general we need the following generalization of a result of Kitaoka:

**Lemma 1.1** *Let $V$ be a vector space over $\mathbf{Q}$ of dimension $m$ with a positive definite quadratic form $q$ on it, let $L_1, \ldots, L_r$ be pairwise non-isometric $\mathbf{Z}$-lattices of rank $m$ on $V$ having the same discriminant. Let $(W_\rho, \rho)$ be a finite dimensional irreducible rational representation of $GL_n(\mathbf{C})$ and let $0 \neq P_i \in \mathcal{H}_q(\rho)$ for $i = 1, \ldots, r$ be $q$-harmonic forms on $V^{m-1}$ with respect to $\rho$ such that $P_i$ is invariant under the action of the group of units $O(L_i)$ of $L_i$.*
*Then the theta series*

$$\vartheta^{(m-1)}(L_i, q, P_i, Z) := \sum_{\mathbf{x} \in L_i^{m-1}} P(\mathbf{x}) \exp(2\pi i tr(q(\mathbf{x})Z))$$

*(with $Z$ in the Siegel upper half plane of degree $m - 1$) are linearly independent.*

*Proof.* With his kind permission we give here a proof due to Y. Kitaoka which is simpler and more elegant than our own original proof.

We show first that given $L$ on $V$ and $0 \neq P \in \mathcal{H}_q(\rho)$ as above there is a global characteristic sublattice $M$ of rank $m - 1$ of $L$ in the sense of Theorem 6.4.1 of [20] such that $P(x_1, \ldots, x_{m-1}) \neq 0$ for a basis (and hence for any basis) $x_1, \ldots, x_{m-1}$ of $M$:

Write $P = \sum_{i=1}^{t} \alpha_i Q_i$ where the $\alpha_i$ are linearly independent over $\mathbf{Q}$ and the nonzero polynomial functions $Q_i : V^{m-1} \to W_\rho$ have rational coefficients (with respect to some fixed bases of $V$ and of $W_\rho$) but are not required to be harmonic forms. We can find an $(m-1)$-tuple $z_1, \ldots, z_{m-1}$ of linearly independent vectors in $L$ such that $Q_1(z_1, \ldots, z_{m-1}) \neq 0$ and put $M' = \mathbf{Z}z_1 + \cdots + \mathbf{Z}z_{m-1}$. Let $S$ be a finite set of primes containing 2 and the primes dividing the discriminant of $L$ and let $\ell$ be some odd prime not dividing the discriminant of $M$. Then by Theorem 6.2.1 of [20] one can, in the same way as in the proof of Proposition 6.4.1 of [20] find linearly independent vectors $x_1, \ldots, x_{m-1}$ in $L$ such that

- The $x_i$ are as close as we want to the $z_i$ in the completion $V_\ell$
- For the $p \in S$ the lattice $\mathbf{Z}_p x_1 + \cdots + \mathbf{Z}_p x_{m-1}$ is a $p$-adic characteristic sublattice of $L_p$ in the sense of [20, p. 153]
- There is some prime $q \notin S \cup \{\ell\}$ such that the discriminant of $M := \mathbf{Z}x_1 + \cdots + \mathbf{Z}x_{m-1}$ is a unit at the primes not in $S \cup \{\ell, q\}$ and exactly divisible by $q$.

The lattice $M$ is then (again as in the proof of Proposition 6.4.1 of [20]) a global characteristic sublattice of $L$, and if we chose the $x_i$ close enough to the $z_i$ at the prime $\ell$, we see that $Q_1(x_1, \ldots, x_{m-1}) \neq 0$ holds, hence also $P(x_1, \ldots, x_{m-1}) \neq 0$, and we have established the existence of the global characteristic sublattice claimed above.

As in the proof of Corollary 6.4.1 of [20] we can now apply the above argument to the $L_i$ to find an index $j$, a characteristic sublattice $M$ of $L_j$ with basis $x_1, \ldots, x_{m-1}$ such that $P_j(x_1, \ldots, x_{m-1}) \neq 0$ and such that $M$ is not represented by any of the $L_i$ with $i \neq j$. By Theorem 6.4.1 of [20] every embedding of $M$ into $L_j$ can be uniquely extended to an isometry in $O(L_j)$. Since $P_j$ was assumed to be $O(L_j)$-invariant, the coefficient at $q((x_1, \ldots, x_{m-1}))$ of $\vartheta^{(m-1)}(L_j, P_j, q, Z)$ is nonzero, whereas that of the other $\vartheta^{(m-1)}(L_i, P_i, q, Z)$ is zero. $\vartheta^{(m-1)}(L_j, P_j, q, Z)$ can therefore not appear in any linear relation between the $\vartheta^{(m-1)}(L_i, P_i, q, Z)$, and iterating this argument we get the assertion.

**Theorem 1.2** *Let $\nu_1 \geq \nu_2 \geq 0$ be integers and let $\tau_1 = \tau_{\nu_1}, \tau_2 = \tau_{\nu_2}$ be the irreducible representations of $D_\infty^\times/\mathbf{R}^\times$ of highest weights $\nu_1, \nu_2$ respectively. Let $\varphi_i \in \mathcal{A}(D_\mathbf{A}^\times, R_\mathbf{A}^\times, \tau_i)$ for $i = 1, 2$ be eigenfunctions of the Hecke operators $\tilde{T}_p$ and of all the involutions $\tilde{w}_p$; assume that $\varphi_1$ and $\varphi_2$ have the same eigenvalues for all the $\tilde{w}_p$.*
*Then if $\nu_1 \neq \nu_2$ and $*$ stands for s the n-th Yoshida-lifting $Y^{(n)}(\varphi_1, \varphi_2, *)$ is defined for $n \geq 2$; if $\nu_1 = \nu_2$ and $*$ stands for $+$ as above, it is defined for $n \geq 1$. the lifting $Y^{(2)}(\varphi_1, \varphi_2, *)$ is in these cases a vector valued Siegel modular form with respect to the representation of highest weight $(\nu_1 + \nu_2, \nu_1 - \nu_2)$ of $GL_2(\mathbf{C})$ and is cuspidal unless*

$\nu_1 = \nu_2$ and $\varphi_1$ and $\varphi_2$ are proportional.
Moreover, under the same conditions $Y^{(3)}(\varphi_1, \varphi_2, *)$ is a nonzero vector valued Siegel modular form for the group $\Gamma_0^{(3)}(N)$ with respect to the representation of highest weight $(\nu_1 + \nu_2 + 2, \nu_1 - \nu_2 + 2, 2)$ of $GL_3(\mathbf{C})$.
If $\nu_1 = \nu_2 = \nu \neq 0$ and $\varphi_1, \varphi_2$ are not proportional, then $Y^{(n)}(\varphi_1, \varphi_2, -)$ is defined for $n \geq 3$ and $Y^{(3)}(\varphi_1, \varphi_2, -)$ is a nonzero vector valued Siegel cusp form for the group $\Gamma_0^{(3)}(N)$ with respect to the representation of highest weight $(2\nu, 3, 3)$ of $GL_3(\mathbf{C})$.
If $\nu_1 = \nu_2 = 0$, then if $\varphi_1, \varphi_2$ are not proportional, $Y^{(3)}(\varphi_1, \varphi_2, -)$ is not defined and $Y^{(4)}(\varphi_1, \varphi_2, -)$ is a nonzero Siegel cusp form for the group $\Gamma_0^{(4)}(N)$ of weight 3.

*Proof.* The claims about the existence of the Yoshida-lifting follow from the description of the theta correspondence at the infinite place in [19]. Lemma 3.1 of [6] gives then the nonvanishing of the function $(\varphi_1 \otimes \varphi_2)_*$ on the adelic orthogonal group $H_\mathbf{A}$. By Lemma 1.2 of [6] the Yoshida-lifting of $(\varphi_1 \otimes \varphi_2)_*$ is a linear combination of theta series of the lattices $I_{ij}$ with harmonic ploynomials that are invariant under the isometries of the lattices $I_{ij}$ equipped with the norm form. The nonvanishing of the Yoshida-lift follows then from Lemma 1.2 above in the first two cases and is obvious in the last case (notice that the nonvanishing of $(\varphi_1 \otimes \varphi_2)_*$ implies in particular that at least one of the $I_{ij}$ has no isometries of determinant $-1$ in this last case).
For the cuspidality properties asserted we notice that from Lemma 1.3 of [6] one sees that Siegel's $\Phi$-operator sends $F := Y^{(n)}(\varphi_1, \varphi_2, *)$ to $Y^{(n-1)}(\varphi_1, \varphi_2, *)$ whenever the latter is defined. From [29] one sees that the $\Phi$-operator sends $F$ to zero if the last two entries in the vector $(\lambda_1, \ldots, \lambda_n)$ describing the highest weight of the representation of $GL_n(\mathbf{C})$ according to which $F$ transforms are not equal. This covers all cases except that $\nu_1 = \nu_2 = \nu$ and $* = -$. If here $\nu = 0$ then $Y^{(4)}(\varphi_1, \varphi_2)$ is scalar valued of weight 3 for $\Gamma_0^{(4)}(N)$, hence annihilated by $\Phi$. For $\nu \neq 0$ and $n = 3$ the form $F|\Phi$ transforms under $\Gamma_0^{(2)}(N)$ according to the representation of highest weight $(2\nu + 2, 3)$, hence is again zero.
The same arguments as in section 8 of [5] apply in the present vector valued situation to show that $Y^{(n)}(\varphi_1, \varphi_2, *)$ is cuspidal if and only if its image under the $\Phi$-operator is zero.

Theorem 1.2 settles the question of the nonvanishing of the Yoshida-lift in the cases where $* = -$. In the remaining cases the argument in the end of the proof reduces the question of the nonvanishing of $Y^{(2)}(\varphi_1, \varphi_2, *)$ to the question of the cuspidality of the nonzero Siegel modular form $Y^{(3)}(\varphi_1, \varphi_2, *)$ in the same way as in [5] for the case $\nu_1 = \nu_2 = 0$. We notice this fact as a corollary:

**Corollary 1.3** *Let the notations be as in Theorem 1.2 and let $*$ denote $+$ if $\nu_1 = \nu_2$ and s otherwise. Moreover $Y^{(2)}(\varphi_1, \varphi_2, *) = 0$ if and only if $Y^{(3)}(\varphi_1, \varphi_2, *)$ is cuspidal.*

2. Vector valued Siegel modular forms

In this section we collect some facts about vector-valued Siegel modular forms; unfortunately there is no adequate reference for the things we need (the notes of Godement [12] are still the best reference for basic facts). In most cases we give only sketches of proof, in particular we omit those calculations which are completely analogous to the scalar-valued case.

**2.1 Differential operators**

**Proposition 2.1** *Let $\rho_0 : Gl(n, \mathbf{C}) \to V_\rho$ be an irreducible polynomial representation with highest weight $(\lambda_1, ..., \lambda_n)$. Then there is a nonzero $V_\rho \otimes V_\rho$-valued polynomial with rational coefficients $Q(s, \mathbf{X}) = Q_s(\mathbf{X})$ in the variables $s$ and $\mathbf{X}$, where $\mathbf{X}$ is a symmetric $2n \times 2n$-matrix of variables $\mathbf{X} = (x_{ij})$, homogeneous of degree $\sum \lambda_i$ as a polynomial in $\mathbf{X}$ with the following property: The associated differential operator $\Delta_s$, given by*

$$(\Delta_s F)(z, w) := (Q_s((\partial_{ij}))F)\left(\begin{pmatrix} z & 0 \\ 0 & w \end{pmatrix}\right)$$

*with $z, w \in \mathbf{H}_n$ maps holomorphic scalar-valued functions $F$ on $\mathbf{H}_{2n}$ to $V_\rho \otimes V_\rho$-valued holomorphic functions on $\mathbf{H}_n \times \mathbf{H}_n$ and satisfies*

$$(2.1) \qquad (\Delta_s) \mid_{\det^s \otimes \rho_0}^{z} M_1 \mid_{\det^s \otimes \rho_0}^{w} M_2 = \Delta_s(F \mid_s M_1^\uparrow M_2^\downarrow)$$

*for all $M_1, M_2 \in Sp(n, R)$, where $M^\uparrow$ and $M^\downarrow$ denote the standard embeddings*

$$\begin{cases} Sp(n, \mathbf{R}) \times Sp(n, \mathbf{R}) & \hookrightarrow & Sp(2n, \mathbf{R}) \\ \begin{pmatrix} A & B \\ C & D \end{pmatrix} \times \begin{pmatrix} a & b \\ c & d \end{pmatrix} & \longmapsto & \begin{pmatrix} A & \mathbf{0} & B & \mathbf{0} \\ \mathbf{0} & a & \mathbf{0} & b \\ C & \mathbf{0} & D & \mathbf{0} \\ \mathbf{0} & c & \mathbf{0} & d \end{pmatrix} \end{cases}$$

*Proof.* It is sufficient to check (2.1) for generators of $Sp(n, \mathbf{R})$ and for "test functions" of type

$$f_{\mathcal{T}} : \mathcal{Z} \longmapsto e^{tr(\mathcal{T} \cdot \mathcal{Z})}$$

with $\mathcal{T} \in \mathrm{Sym}_{2n}(\mathbf{C})$ and $\mathcal{Z} = \begin{pmatrix} z & z_2 \\ z_2' & w \end{pmatrix} \in \mathbf{H}_{2n}$ We fix the "weight" $s$ and consider $Q_s$ as a polynomial in the $\partial_{ij}$ with unknown coefficients (depending on $s$). For generators of type $\begin{pmatrix} 1_n & R \\ 0_n & 1_n \end{pmatrix}$, $R \in \mathrm{Sym}_n(\mathbf{R})$ the equation (2.1) is satisfied, so we only have to investigate $I^\uparrow$, $I^\downarrow$ for $I = \begin{pmatrix} 0_n & -1_n \\ 1_n & 0_n \end{pmatrix}$. We have to compare

$$(2.2) \qquad \Delta_s(f_{\mathcal{T}} \mid_s I^\uparrow)$$

and

$$(2.3) \qquad \det(z)^{-s} \cdot (\rho_0(z) \otimes id_{V_{\rho_0}})^{-1} \cdot (\Delta_s f_{\mathcal{T}})(I^\uparrow <\mathcal{Z}>)$$

(and a similar expression with $I^\downarrow$ instead of $I^\uparrow$). After multiplication by $\rho_0(z) \otimes id_{V_{\rho_0}}$ and a suitable power of $\det(z)$ both (2.2) and (2.3) are of the form

$$e^{tr(-T_1 \cdot z^{-1} + T_4 \cdot w)} \times (\text{polynomial in } z \text{ and } \mathcal{T})$$

Here of course $T_1, ... T_4$ denote the blocks of size $n$ in the matrix T. From this we see that (2.1) for all $f_{\mathcal{T}}$ and $(M_1, M_2) = (I, id)$ or $(id, I)$ is equivalent to a system $\mathcal{L}_s$ of linear equations for the coefficients of $Q_s$. The coefficients of this system depend polynomially on $s$. Ibukiyama [17] proved that the space of solutions of $\mathcal{L}_s$ is one-dimensional for $s = k \in \mathbf{N}, k \geq n$. Therefore the corank of $\mathcal{L}_s$ is generically in $s$ equal to 1. The generic solution of $\mathcal{L}_s$ (properly normalized, such that the coefficients of $\partial_{ij}$ are polynomials in $s$ with no factor in common) gives us the existence of $\Delta_s$; $\Delta_s$ becomes unique, if we choose some monomial in the $\partial_{ij}$ with non-zero coefficient and require this coefficient to be a monic polynomial in $s$.

Remark: The polynomial $Q_s$ automatically satisfies the symmetry relation

$$(2.4) \qquad Q_s(\partial_{ij}) = Q_s\left((\partial_{ij}) \left[\begin{pmatrix} \mathbf{0}_n & \mathbf{1}_n \\ \mathbf{1}_n & \mathbf{0}_n \end{pmatrix}\right]\right)$$

thanks to the result of Ibukiyama.

On $V_{\rho_0}$ we have a non-degeneate bilinear form $<,>$ such that

$$(2.5) \qquad <v_1, \rho(g)v_2> = <\rho(g^t)v_1, v_2>$$

We denote by $\bullet : V_\rho \otimes V_\rho \longmapsto V_\rho \otimes V_\rho^\star = \text{End}(V_\rho)$ the isomorphism given by

$$(2.6) \qquad (v_1 \otimes v_2)^\bullet(v) = <v, v_1> v_2$$

Then we have for all $g_1, g_2 \in Gl(n, \mathbf{C})$

$$(2.7) \qquad (\rho_o(g_1)v_1 \otimes \rho_o(g_2)v_2)^\bullet = \rho_0(g_1)(v_1 \otimes v_2)^\bullet \rho_o(g_2^t)$$

Denoting by $\mathbf{h}_s$ the function

$$h_s : \begin{cases} \mathbf{H}_n \longrightarrow \mathbf{C} \\ \mathcal{Z} \longmapsto \det(z_1 + z_2 + z_2^t + z_4)^{-s} \end{cases}$$

we claim the following important relation

$$(2.8) \qquad \Delta_s^\bullet \mathbf{h}_s = c(s) \cdot \rho_0(z_1 + z_4)^{-1} \cdot \det(z_1 + z_4)^{-s}$$

with a certain scalar $c(s) = c_{\rho_0}(s)$.

The *proof* of this claim is quite similar to the proof of Théoréme 5 in [12, Exp.6]: The function $\mathbf{h}_s$ is rich of symmetries, which thanks to (2.1) gives rise to enough symmetry properties of $\Delta_s^\bullet$ to prove (2.8). By considering $\Delta_s^\bullet \mathbf{h}_s$ as a function $\mathbf{f}_s(z_1 + z_4, z_1 - z_4)$ of the variables $z_1 + z_4$ and $z_1 - z_4$ we see (by applying suitable translations, using (2.1))

$$\mathbf{f}_s(z_1 + z_4, z_1 - z_4 + 2T) = \mathbf{f}_s(z_1 + z_4, z_1 - z_4)$$

for all $T \in \mathrm{Sym}_n(\mathbf{R})$, hence

$$\Delta_s^\bullet \mathbf{h}_s = \mathbf{g}_s(z_1 + z_4)$$

with a suitable function $\mathbf{g}$ on $\mathbf{H}_n$. For $g \in Gl(n, \mathbf{R})$

$$\mathbf{h}_s \mid_s \begin{pmatrix} g^t & 0 \\ 0 & g^{-1} \end{pmatrix}^\uparrow \mid_s \begin{pmatrix} g^t & 0 \\ 0 & g^{-1} \end{pmatrix}^\downarrow = \mathbf{h}_s$$

hence

$$\Delta_s^\bullet(\mathbf{h}_s)(z_1, z_4) = \det(g)^{2s} \rho_0(g)(\Delta_s^\bullet \mathbf{h}_s)(g^t z_1 g, g^t z_4 g) \cdot \rho_0(g^t)$$

in particular

$$\mathbf{g}_s(iy) = \det(y)^{-s} \cdot \rho_0(y^{-\frac{1}{2}}) \cdot \mathbf{g}_s(i\mathbf{1}_n) \cdot \rho_0(y^{-\frac{1}{2}})$$

It remains to show that $\mathbf{g}_s(i\mathbf{1})$ is of the form $scalar \times id_{V_{\rho_o}}$. We start from a symmetry relation for $\mathbf{h}_s$, namely

(2.9) $$\mathbf{h}_s \mid_s M^\uparrow = \mathbf{h}_s \mid_s \tilde{M}^\downarrow$$

for all $M \in Sp(n, \mathbf{R})$, where $\tilde{M} = M^{-1} \begin{bmatrix} \mathbf{1}_n & \mathbf{0}_n \\ \mathbf{0}_n & -\mathbf{1}_n \end{bmatrix}$ This relation (2.9) is again proved by considering generators of $Sp(n, \mathbf{R})$, the of translations being trivial. For M = I we prove

(2.10) $$\mathbf{h}_s \mid_s (I^\uparrow \circ (-I)^\downarrow) = \mathbf{h}_s$$

To prove (2.10) we use the (elementary) relation

$$\mathbf{h}_s = \mathcal{G}_s \mid_s I^\uparrow = \mathcal{G}_s \mid_s I^\downarrow$$

with

$$\mathcal{G}_s(\mathcal{Z}) = \det\left(\begin{pmatrix} \mathbf{0}_n & \mathbf{1}_n \\ \mathbf{1}_n & \mathbf{0}_n \end{pmatrix} \cdot \mathcal{Z} + \begin{pmatrix} \mathbf{1}_n & \mathbf{0}_n \\ \mathbf{0}_n & \mathbf{1}_n \end{pmatrix}\right)^{-s}$$

We obtain

$$\begin{aligned}
\mathbf{h}_s \mid_s (I^\uparrow \circ (-I)^\downarrow) &= \mathcal{G}_s \mid_s \left(I^\uparrow \circ I^\uparrow \circ (-I)^\downarrow\right) \\
&= \mathcal{G}_s \mid_s \left((-\mathbf{1}_{2n})^\uparrow \circ (-I)^\downarrow\right) \\
&= \mathcal{G}_s \mid \left(\left(I^\downarrow\right) \circ (-\mathbf{1}_{4n})\right) \\
&= \mathcal{G}_s \mid_s I^\downarrow \\
&= \mathbf{h}_s
\end{aligned}$$

Therefore we have for any $M = \begin{pmatrix} a & b \\ c & d \end{pmatrix} \in Sp(n, \mathbf{R})$

$$\Delta_s(\mathbf{h}_s \mid_s M^\uparrow \circ \tilde{M}^{-1\downarrow}) = \Delta_s(\mathbf{h}_s)$$

which implies

(2.11) $(\det(cz_1 + d)\det(-cz_4 + d))^{-s}$
$$\times \rho_0\left((cz_1 + d)^{-1}\right)\mathbf{g}_s\left(M<z_1> + \tilde{M}^{-1}<z_4>\right)\rho_0\left((-cz_4 + d)^t\right)^{-1} = \mathbf{g}_s(z_1 + z_4)$$

We apply (2.11) to $M = \begin{pmatrix} u & -v \\ v & u \end{pmatrix} \in Sp(n, \mathbf{R})$ with $u + iv \in U(n, \mathbf{C})$ and $z_1 = z_4 = i \cdot \mathbf{1}_n$:

(2.12) $$\rho_0(u+iv)^{-1}\mathbf{g}_s(2i)\rho_0\left((u-iv))^t\right)^{-1} = \mathbf{g}_s(2i)$$

This means that $\mathbf{g}_s(2i)$ must be a scalar matrix, because it commutes with all $\rho_0(g)$, $g \in U(n, \mathbf{C})$; this follows from (2.12) because of

$$(u - iv)^t = (u + iv)^{-1}$$

The proof of (2.8) is finished.

### 2.2 Pullbacks of Eisenstein Series

Using the differential operator $\Delta_s$ we extend the machinery of "pullbacks of Eisenstein series" to vector valued modular forms. The scalar valued case was previously discussed in [11, 1, 5] and the case of symmetric tensor representations in [3] (and recently in [28]). Our approach is somewhat different fromTakayanagi's [28], who considers a differential operator independent of s; we avoid some of the combinatorial problems coming up in [28]. We denote by $\{,\}$ a hermitian scalar product on $V_\rho$ such that $\{\rho(\bar{g}^t)v_1, v_2\} = \{v_1, \rho(g)v_2\}$; after a suitable normalization we may assume

(2.13) $$\{v_1, v_2\} = <v_1, v_2^\sigma>$$

where $\sigma$ denotes the natural complex conjugation on $V_\rho$.
For $F \in [\Gamma_0^n(N), \rho]_0$ with $\rho = \rho_0 \otimes \det^k$ we consider

(2.14) $$\int_{\Gamma_0(N)\backslash \mathbf{H}_n} \left\{\rho(\sqrt{v})F(w), \rho(\sqrt{v})\Delta_{\bar{s}+k}\mathbf{G}_{2n}^k\left(\begin{pmatrix} w & 0 \\ 0 & -\bar{z} \end{pmatrix}, \bar{s}\right)\right\} \det(v)^s \det(y)^s d\omega_w$$

where
$$\mathbf{G}_{2n}^k(\mathcal{Z}, s) = \sum_{M \in \Gamma_\infty^{2n} \backslash \Gamma_0^{2n}(N)} j((M, \mathcal{Z})^{-k-s} j(M, \bar{\mathcal{Z}})^{-s}$$

(2.15)
$$= \sum_{M = \begin{pmatrix} \star & \star \\ C & D \end{pmatrix}} \det(C\mathcal{Z} + D)^{-k-s} \det(C\bar{\mathcal{Z}} + D)^{-s}$$

To compute this integral, one follows the same line of computation as in [5, section 1], using the double coset decomposition for
$$\Gamma_\infty^{2n} \backslash \Gamma_0^{2n}(N) / \Gamma_0^n(N)^\uparrow \times \Gamma_0^n(N)^\downarrow$$
as described in [5, Thm.1.1], as well as the description of the corresponding left coset decomposition [5, Thm.1.2]. The standard unfolding procedure yields (using (2.1))

(2.16) $$\sum_M \left( \int_{\mathbf{H}_n} \left\{ \rho(v) F(w), (\Delta_{\bar{s}+k} \mathbf{h}_{\bar{s}+k})(w, -\bar{z}) \det(\bar{w} - z)^{-\bar{s}} \right\} \det(v)^s \det(y)^s d\omega \right) |T_N(M) \det(M)^{-k-s}$$

where the summation is over all integral diagonal matrices
$$M = \begin{pmatrix} m_1 & & \\ & \ddots & \\ & & m_n \end{pmatrix}$$
with $M \equiv 0$ modulo $N$ and $0 < m_1 | m_2 \dots | m_n$, and where $T_N(M)$ denotes the Hecke operator associated to the double coset
$$\Gamma_0(N) \begin{pmatrix} \mathbf{0}_n & -M^{-1} \\ M & \mathbf{0}_n \end{pmatrix} \Gamma_0(N).$$

In view of (2.13) and (2.8) it remains to investigate
(2.17)
$$\int_{\mathbf{H}_n} \det(-\bar{z} + w)^{-s} ((\Delta_{\bar{s}+k} \mathbf{h}_{\bar{s}+k}(w, -\bar{z}))^\sigma)^\bullet (\rho(v) F(w)) \det(v)^s \det(y)^s d\omega_w$$
$$= \overline{c(\bar{s}+k)} \int_{\mathbf{H}_n} \overline{\rho_0(-\bar{z}+w)^{-1} \det(-\bar{z}+w)^{-k-\bar{s}}} \det(-\bar{z}+w)^{-s} \rho(v) F(w) \det(v)^s \det(y)^s d\omega_w$$
$$= \overline{c(\bar{s}+k)} \int_{\mathbf{H}_n} \rho(-z+\bar{w})^{-1} \det(-z+\bar{w})^{-k-s} \det(-\bar{z}+w)^{-s} \rho(v) F(w) \det(v)^s \det(y)^s d\omega_w$$

From [12, Exp.6,Thm.5], applied to the holomorphic(!) function
$$w \longmapsto \det(-\bar{z}+w)^{-s} F(w)$$

with $Re(s) \gg 0$ it is clear that (2.17) is a scalar multiple of $F(z)$; we denote this scalar by $\mathcal{H}_{\rho_0}(s+k)$.

It is actually possible to compute this scalar explicitly by transferring the integral to the generalized unit disc

$$\mathcal{D}_n = \{\xi \in \mathrm{Sym}_n(\mathbf{C}) | \mathbf{1}_n - \bar{\xi}\xi > 0\}$$

by a Cayley-transform (see e.g.[21]); as a result we get

$$\mathcal{H}_{\rho_0}(s+k) = 2^{n^2+n-4ns}(\frac{1}{2}i)^{nk+\sum \lambda_i} \int_{\mathcal{D}_n} \rho_o \otimes \det(\mathbf{1}_n - \bar{\xi}\xi)^{s+k} \frac{d\xi}{\det(\mathbf{1}_n - \bar{\xi}\xi)^{n+1}}$$

The latter integral (which is easily seen to be of type scalar$\times id_{V_\rho}$) is a generalized version of the integrals computed in [16]; using the methods of [16] one can explicitly compute this integral, but we do not do this here.

What we shall later on need is a relation between $c_{\rho_0}(s) \cdot \mathcal{H}_{\rho_0}(s)$ and the polynomial $Q_s\left(\begin{pmatrix} \mathbf{1}_n & \mathbf{1}_n \\ \mathbf{1}_n & \mathbf{1}_n \end{pmatrix}\right)$:

For $s \in \mathbf{C}$ with $Re(s) \gg 0$ we start from the well-known identity

(2.18)
$$\sum_{L \in \mathrm{Sym}_n(\mathbf{Z})} \det(z_1 + z_2 + z_2^t + z_4)^{-s} = \frac{A}{\Gamma_n(s)} \sum_{T \in \Lambda_n^+} \det(T)^{s-\frac{n+1}{2}} e^{2\phi itr(T(z_1+z_2+z_2^t+z_4))}$$

where $\Lambda_n^+$ denotes the set of all half-integral symmetric positive matrices of size $n$, $\mathcal{Z} = \begin{pmatrix} z_1 & z_2 \\ z_2^t & z_4 \end{pmatrix}$ is an element of $\mathbf{H}_{2n}$ and

$$A = A(s) = e^{-\frac{1}{2}\pi ins} 2^{n(s-\frac{1}{2}(n-1))} \pi^{sn}$$

We apply $\Delta_s^\bullet$ to both sides of (2.18):

(2.19)
$$\frac{A}{\Gamma_n(s)} \sum_{T \in \Lambda_n^+} \det(T)^{s-\frac{n+1}{2}} Q_s\left(\begin{pmatrix} T & T \\ T & T \end{pmatrix}\right)^\bullet e^{2\pi itr(T(z_1+z_4))}$$
$$= c_{\rho_o}(s) \sum_{L \in \mathrm{Sym}_n(\mathbf{C})} \rho_o(z_1 + z_4 + L)^{-1} \det(z_1 + z_4 + L)^{-s}$$

According to [12, Exp.10,Thm.7,Corollaire] this is equal to

$$c_{\rho_0}(s)\mathcal{H}_{\rho_0}(s) \sum_{T \in \Lambda_n^+} H_{\rho_0 \otimes \det^s}(4T)^{-1} e^{2\pi itr(T(z_1+z_4))}$$

with

(2.20)
$$H_{\rho_0 \otimes \det^s}(T) = \int_{Y^{(n)}>0} \rho_0(Y) \det(Y)^s e^{-2\pi tr(TY)} \det(Y)^{-n-1} dY$$

Therefore we get for all $T \in \Lambda_n^+$ and $Re(s) \gg 0$ the relation

$$(2.21) \quad \frac{A}{\Gamma_n(s)} \det(T)^{s-\frac{n+1}{2}} Q_s \left( \begin{pmatrix} T & T \\ T & T \end{pmatrix} \right)^\bullet \circ H_{\rho_o \otimes \det^s}(4T) = c_{\rho_o}(s) \mathcal{H}_{\rho_o}(s) \cdot id_{V_{\rho_o}}$$

It is obvious that $H_{\rho_0 \otimes \det^s}(T)$ has a meromorphic continuation to all of $\mathbf{C}$ since it can be obtained from the function

$$T \longmapsto \int_{Y^{(n)} > 0} \det(Y)^{s-n-1} e^{-2\pi i tr(TY)} dY = \Gamma_n(s - \frac{n+1}{2}) \det(2\pi T)^{-s + \frac{n+1}{2}}$$

by applying a suitable differential operator (w.r.t.$T$). In particular, equation (2.21) provides the meromorphic continuation of $c_{\rho_0}(s) \cdot \mathcal{H}_{\rho_0}(s)$ to all of $\mathbf{C}$ (if we do not prefer to compute $\mathcal{H}_{\rho_0}$ explicitly).

Later on we need some information about poles of $H_{\rho_0 \otimes \det^s}(S)$, which we provide here ($S = S^t > 0$): Let $v_0 \in V_{\rho_0}$ be a vector of highest weight, then

$$<v_0, H_{\rho_0 \otimes \det^s}(S) v_0> = \int_{Y>0} <v_0, \rho_0(Y) v_0> \det(Y)^{s-(n+1)} e^{-\pi tr(SY)} dY$$

$$= \int_{Y>0} \left( \prod_{i=1}^n t_{ii}^{2\lambda_i} \right) \det(Y)^{s-(n+1)} e^{-\pi tr(SY)} dY$$

with $Y = T^t \cdot T$, $T$ upper triangular. This kind of integrals is studied in [22, p.76], one obtains

(2.22)
$$<v_0, H_{\rho_0 \otimes \det^s} v_0> = \pi^{\frac{1}{4}n(n-1)} \left( \prod_{i=1}^n \Gamma(\lambda_i + s - \frac{i}{2} - \frac{n}{2}) \right) \cdot \left( \prod_{i=1}^n r_{ii}^{\lambda_i + s - \frac{n+1}{4} - \frac{i}{2}} \right)$$

where $(\pi S)^{-1} = R \cdot R^t$ with R upper triangular.

**2.3 Hecke Eigenforms**

We assume now that $F \in [\Gamma_0^n(N), \rho_0 \otimes \det^k]_0$ is an eigenform of all Hecke operators $T_N(M)$, the corresponding eigenvalue will be denoted by $\lambda_F(M)$. If the eigenvalue $\lambda_F(N \cdot \mathbf{1}_n)$ is different from zero, we get the following Euler product (see [5, Thm.2.1])

$$\sum_M \lambda_F(M) \det(M)^{-s} = \frac{\lambda_F(N \cdot \mathbf{1}_n)}{N^{ns}} \cdot \frac{1}{\zeta^{(N)}(s) \prod_{i=1}^n \zeta^{(N)}(2s - 2i)} \cdot \Lambda_N(s-n) \cdot D_F^{(N)}(s-n)$$

where

$$D_F^{(N)}(s) = \prod_{p \nmid N} \left( \frac{1}{1 - p^{-s}} \cdot \prod_{i=1}^n \frac{1}{(1 - \alpha_{ip} p^{-s})(1 - \alpha_{ip}^{-1} p^{-s})} \right)$$

denotes the (restricted) standard L-function attached to F, the $\alpha_{ip}$ being the Satake parameters attached to F in the usual way; moreover

$$\Lambda_N(S) = \prod_{q|N} \prod_{i=1}^{n} \frac{1}{1 - \beta_{iq} q^{-s}}$$

with $\beta_{iq}$ being the Satake parameters describing the following one-dimensional representation of the Hecke-algebra attached to the Hecke pair $\left(Gl(n, \mathbf{Z}), Gl(n, \mathbf{Z}[\frac{1}{q}])\right)$:

$$Gl(n, \mathbf{Z}) M_0 Gl(n, \mathbf{Z}) \longmapsto \frac{1}{\lambda_F(N \cdot \mathbf{1}_n)} \cdot \lambda_F(N \cdot M_0)$$

for $M_0 \in \mathbf{Z}^{n,n)}$, $\det(M_0) =$ power of $q$. For details we refer to [5, section 2]. Summarizing our computation, we obtain for such a Hecke eigenform $F \in [\Gamma_0^n(N), \rho]_0$, $\rho = \rho_0 \otimes \det^k$:

$$(2.23) \quad \int_{\Gamma_0(N)\backslash \mathbf{H}_n} \left\{ \rho(\sqrt{v}) F(w), \rho(\sqrt{v}) \Delta_{\bar{s}+k} \mathbf{G}_{2n}^k \left( \begin{pmatrix} w & 0 \\ 0 & -\bar{z} \end{pmatrix}, \bar{s} \right) \right\} \det(v)^s \det(y)^s d\omega^w$$

$$= \overline{c_{\rho_0}(\bar{s}+k)} \cdot \mathcal{H}_{\rho_0}(s+k) \cdot \frac{\lambda_F(N \cdot \mathbf{1}_n)}{(N^n)^{k+2s}} \cdot \frac{\Lambda_N(2s+k-n) \cdot D_F^{(N)}(2s+k-n)}{\zeta^{(N)}(2s+k) \prod_{i=1}^{n} \zeta^{(N)}(4s+2k-2i)}$$

## 2.4 The scalar-valued case revisited

The scalar-valued case ($\rho_0 = \det^\nu$) can be made more explicit; the differential operators $\Delta_s = \Delta^{(\nu)}$ in this case were constructed in a different way in [2]. We normalize them by fixing the coefficient of $\det(\partial_2)$ with $(\partial_{1 \le i,j \le 2n}) = \begin{pmatrix} \partial_1 & \partial_2 \\ \partial_2^t & \partial_4 \end{pmatrix}$ to be

$$C_n(s - n + \frac{1}{2} + \nu - 1) \ldots C_n(s - n + \frac{1}{2})$$

with

$$C_n(s) = s(s + \frac{1}{2}) \ldots (s + \frac{n-1}{2}) = \frac{\Gamma_n(s + \frac{n+1}{2})}{\Gamma_n(s + \frac{n-1}{2})}$$

With this normalization,

$$(2.24) \quad c_\rho(s) = c_\nu(s) = \frac{\Gamma_n(s + \nu) \cdots + \nu - \frac{n}{2})}{\Gamma_n(s) \cdot \Gamma_n(s - \frac{n}{2})}$$

Moreover
(2.25)
$$\begin{aligned}\mathcal{H}_{\rho_o}(s+k) &= \mathcal{H}_\nu(s+k) \\ &= 2^{n^2+n-4ns-nk-n\nu} \cdot i^{n(n+\nu)} \int_{D_n} \det(1-\xi\bar{\xi})^{s+k+\nu-n-1} d\xi \\ &= i^{n(k+\nu)} \cdot 2^{n^2+n-4ns-nk-n\nu-\frac{n(n+1)}{2}} \cdot \frac{\Gamma_n(k+\nu+s-\frac{n+1}{2})}{\Gamma_n(k+\nu+s)}.\end{aligned}$$

We need a slightly modified version of (2.23) later on (still in the scalar-valued case):

In [4] it is described, how one gets (essentially from a variant of Shimura's theory of nearly holomorphic functions) an operator identity (with $p_k(s)$, $d_k(s)$ being polynomials in $s$)

(2.26) $$p_k(s) \cdot \det(y_1)^s \det(y_4)^s \cdot \Delta^\nu_{s+k} \left( \det(Y)^{-s} \times ? \right) = d_k(s) \cdot \Delta^\nu_k + \mathcal{R}$$

Here $\mathcal{R}$ is a finite sum of differential operators of type

(2.27) $$(\mathbf{D}^{z_1} \otimes \mathbf{D}^{z_4}) \circ \mathcal{D}$$

where $\mathcal{D}$ is some holomorphic differential operator of the type considered in section 2.1 and $\mathbf{D}^{z_1}$ ($\mathbf{D}^{z_4}$ respectively) is a Maaß type differential operator with respect to $z_1 \in \mathbf{H}_n$ ($z_4 \in \mathbf{H}_n$ respectively); in (2.27) at least one of the Maaß type operators is non-trivial. From [4] we quote the formula

(2.28) $$\frac{d_k(s)}{p_k(s)} = \prod_{\mu=1}^{\nu} \frac{C_n(k+s-n+\nu-\frac{\mu}{2})}{C_n(k-n+\nu-\frac{\mu}{2})}$$

Moreover we use that holomorphic Siegel modular forms are orthogonal (with respect to the Petersson scalar product) to automorphic forms in the image of the Maaß operators. For details on these facts we refer to [4]. For an eigenform $F \in \left[\Gamma_o(N), \det^{k+\nu}\right]_0$ we obtain now

(2.29)
$$\int_{\Gamma_0(N)\backslash\mathbf{H}_n} \int_{\Gamma_0(N)\backslash\mathbf{H}_n} F(w) \overline{\Delta_k \left( E^k_{2n}\left( \begin{pmatrix} w & \mathbf{0} \\ \mathbf{0} & -\bar{z} \end{pmatrix}, \bar{s} \right) \right) \det(v)^{k+\nu} d\omega_w} F(z) \det(y)^{k+\nu} d\omega_z$$

$$= \overline{c_\nu(\bar{s}+k)} \cdot \mathcal{H}_\nu(s+k) \cdot \frac{d_k(s)}{p_k(s)} \cdot \frac{\lambda_F(N \cdot \mathbf{1}_n)}{(N^n)^{k+2s}} \cdot \frac{\Lambda_N(2s+k-n) \cdot D^{(N)}_F(2s+k-n)}{\zeta^{(N)}(2s+k) \cdot \prod_{i=1}^{n} \zeta^{(N)}(4s+2k-2i)} \cdot \langle F, F \rangle.$$

Here of course
$$\begin{aligned}E^k_{2n}(\mathcal{Z}, s) &= \det(Y)^s \cdot G^k_{2n}(\mathcal{Z}, s) \\ &= \sum_{M \in \Gamma^{2n}_\infty \backslash \Gamma^{2n}_0(N)} \det(Y)^s \mid_k M\end{aligned}$$

3. The key proposition

The key to non-vanishing properties of Yoshida-liftings is

**Proposition 3.1** *Let $0 \neq F \in \left[\Gamma_0^3(N), \rho_0 \otimes \det^2\right]_0$ be an eigenform of all $T_N(M)$, $M \equiv 0 \mod N$ and assume that $\rho_0$ is of highest weight $(\lambda_1, \lambda_2, 0)$. Then $\Lambda_N(s) \cdot D_F^{(N)}(s)$ is regular in $s = 1$.*

*Proof:* As in [5, section 9] we start from the fact that $G_6^2$ has a simple pole in $s = 1$. In view of (2.23) it is sufficient to show that $c_{\rho_0}(s) \cdot \mathcal{H}_{\rho_0}(s)$ has a pole (of any order) in $s = 3$; using (2.21) we investigate $Q_s^\bullet\left(\begin{pmatrix} T & T \\ T & T \end{pmatrix}\right) \cdot H_{\rho_0 \otimes \det^s}(4T)$ instead (for any $T \in \Lambda_n^+$):

From (2.22) we see that $H_{\rho_0 \otimes \det^s}(4T)$ must have a pole in $s = 3$ (i.e. in some matrix realization of $H_{\rho_0 \otimes \det^s}(4T)$ there is at least one coefficient in this matrix having a pole in $s = 3$. We are done, if we can show that $Q_s^\bullet\left(\begin{pmatrix} T & T \\ T & T \end{pmatrix}\right)$ has determinant different from zero in $s = 3$, when considered as an element of $End(V_{\rho_0})$. It is notationally easier to work with the (real) vectorspace

$$V_{\rho_0}^+ = \{v \in V_{\rho_0} \mid v = v^\sigma\}$$

for the moment. The statement in question about $\det(Q_s^\bullet)$ is true more generally:

**Propositon 3.2** *Let $l, n$ be natural numbers with $l \geq 2n$, $\rho_0$ a polynomial representation of $Gl(n, \mathbf{C})$ and consider the polynomial function*

$$\mathbf{P}_l : \begin{cases} \mathbf{R}^{(l,n)} \times \mathbf{R}^{(l,n)} & \longrightarrow \quad V_{\rho_0}^+ \otimes V_{\rho_0}^+ \\ (X, Y) & \longmapsto \quad Q_l(X^t \cdot X, X^t \cdot Y, Y^t \cdot Y) \end{cases}$$

*Then there is a $X_0 \in \mathbf{R}^{(l,n)}$ such that $\det(\mathbf{P}_l(X_0, X_0)^\bullet) \neq 0$.*

It is clear (by continuity) that $X_0$ may be chosen to be of maximal rank; by the action of $Gl(n, \mathbf{R})$ from the left and the right on $Q_l$ this implies that $Q_l\left(\begin{pmatrix} T & T \\ T & T \end{pmatrix}\right)$ has non-zero determinant for **all** positive definite $T \in \text{Sym}_n(\mathbf{R})$.

*Proof.* As a polynomial in $X$ and $Y$, $\mathbf{P}_l$ is an element of $Harm(l, n) \otimes_{\mathbf{C}} Harm(l, n)$ and also of $Harm^+(l, n) \otimes_{\mathbf{R}} Harm^+(l, n)$; by its construction (see [17] and Prop.2.1), $\mathbf{P}_l$ is different from zero and

(3.1) $$f \longmapsto << f, \mathbf{P}_l(*, Y) >>$$

defines an endomorphism of $Harm_{\rho_0}^+(l, n)$ commuting with the action of $O(n, \mathbf{R})$; here $<<, >>$ denotes the natural scalar product on $Harm_{\rho_0}^+(l, n)$ "adapted" to $O(l, \mathbf{R})$.

The endomorphism (3.1) is selfadjoint (because $\mathbf{P}_l$ is symmetric), hence, by Schur's lemma, it is scalar, i.e. there is a constant $c \neq 0$ such that

$$(3.2) \qquad \mathbf{P}_l(X,Y) = c \cdot \sum_{i=1}^{r} f_i(X) \otimes f_i(Y)$$

where $f_i$ runs through some orthonormal basis of the $r$-dimensional space

$$Harm^+_{\rho_0}(l,n).$$

Now we recall from [10] that $Harm^+_{\rho_0}(l,n)$ is generated by functions of type

$$(3.3) \qquad X \longmapsto \rho_0(X^t \cdot A)v$$

with $v \in V_{\rho_0}$, $A \in \mathbf{C}^{(l,n)}$ with $A^t \cdot A = 0$ Therefore there exists $X_0 \in \mathbf{C}^{(l,n)}$ such that

$$(3.4) \qquad \{f(X_0) \mid f \in Harm_{\rho_0}(l,n)\} = V_{\rho_0}$$

(take $X_0$ in (3.3) such that $X_0^t \cdot A$ has maximal rank).

We can reformulate (3.4) by saying that there is $X_0$ such that the matrix

$$(<f_i(X_0), v_j>_{ij})$$

has maximal rank, where $v_j$ runs through a basis of $V_{\rho_0}$ (or $V^+_{\rho_0}$). From this it follows immediately that there exists a $X_0 \in \mathbf{R}^{(l,n)}$ with

$$(3.5) \qquad \left\{f(X_0) \mid f \in Harm^+_{\rho_0}(l,n)\right\} = V^+_{\rho_o}$$

Now let $X_0$ be as in (3.5); then for $v \in V_{\rho_0}$ the equation $\mathbf{P}_l(X_0,X_0)^\bullet(v) = 0$ implies $\sum_i <f_i(X_0), v> f_i(X_0) = 0$ and hence $\sum_i |<f_i(X_0), v>|^2 = 0$, so $v$ must be zero.

## 4. Hecke action on Yoshida liftings

Throughout this section we let $\nu_1 \geq \nu_2 \geq 0$ as in Section 1 and let $\varphi_i \in \mathcal{A}(D^\times_\mathbf{A}, R^\times_\mathbf{A}, \tau_{\nu_i})$ for $i = 1, 2$ be given. By $\mathcal{A}_{\mathrm{ess},p}(D^\times_\mathbf{A}, R^\times_\mathbf{A}, \tau_{\nu_i})$ we denote the $p$-essential part of $\mathcal{A}(D^\times_\mathbf{A}, R^\times_\mathbf{A}, \tau_{\nu_i})$, i. e. , the set of $\varphi \in \mathcal{A}(D^\times_\mathbf{A}, R^\times_\mathbf{A}, \tau_{\nu_i})$ which are orthogonal to all functions invariant under an order $\tilde{R}$ with $\tilde{R}_p$ strictly containing $R_p$. We assume that for each $p \mid N$ at least one of the $\varphi_i$ belongs to the $p$-essential part of $\mathcal{A}(D^\times_\mathbf{A}, R^\times_\mathbf{A}, \tau_{\nu_i})$. We assume moreover that $\varphi_1, \varphi_2$ are nonproportional eigenfunctions of all the Hecke operators $\tilde{T}(p)$ for the primes $p \nmid N$ with eigenvalues $\lambda^{(i)}_p$ and eigenfunctions of the involutions $\tilde{w}_p$ for the $p \mid N$ with the same eigenvalues. (If instead of these assumptions we start out with two essential eigenforms $\varphi' \in \mathcal{A}(D^\times_\mathbf{A}, (R'_\mathbf{A})^\times, \tau_{\nu_1})$ and $\varphi'' \in \mathcal{A}(D^\times_\mathbf{A}, (R''_\mathbf{A})^\times, \tau_{\nu_2})$ with Eichler orders $R', R''$ of square free levels $N', N''$ in $D$ then we can assume without loss of generality that $R' \cap R''$ is an Eichler order of level $N = \mathrm{lcm}(N', N'')$ and go over from $\varphi', \varphi''$ to forms $\varphi_1, \varphi_2$ that satisfy the assumptions stated above and have the same Hecke eigenvalues as $\varphi', \varphi''$ respectively for the $p \nmid N$.)

Since the question of nonvanishing of the Yoshida-lift $Y^n(\varphi_1, \varphi_2, -)$ in the case $\nu_1 = \nu_2$ has been settled in Section 1, we let $Y^3(\varphi_1, \varphi_2)$ denote $Y^3(\varphi_1, \varphi_2, +)$ or $Y^3(\varphi_1, \varphi_2, s)$, depending on whether $\nu_1 = \nu_2$ holds or not.

According to the results of the previous section we have to show now that the Yoshida-lift $Y^3(\varphi_1, \varphi_2)$ is an eigenfunction of all Hecke-operators for the $p \nmid N$ as well as for the operators $T_N(\mathbf{M})$ for the bad places and to calculate the Euler factors for $p \nmid N$ of the standard $L$-function of the Yoshida-lift as well as the factor $\Lambda_N(s)$ for the bad places. The results and proofs are here completely analogous to those of [5] for the case $\nu_1 = \nu_2 = 0$ (as is clear from an adelic point of view), and we sketch them only very briefly.

**Theorem 4.1** *Let $\varphi_i$ for $i = 1, 2$ be as above with different Hecke eigenvalues and put $k_i = 2 + 2\nu_i$. For primes $p \nmid N$ denote by $\beta_p, \beta_p^{-1}$ respectively $\widetilde{\beta}_p, \widetilde{\beta}_p^{-1}$ the Satake parameters of $\varphi_1$ respectively $\varphi_2$ with respect to the Hecke-algebra of $GL_2(\mathbf{Q}_p) \cong D_p^\times$ (so that $\lambda_p^{(1)} = p^{(k_1-1)/2}(\beta_p+\beta_p^{-1})$, $\lambda_p^{(2)} = p^{(k_2-1)/2}(\widetilde{\beta}_p+\widetilde{\beta}_p^{-1})$). Then for $*$ denoting $+, -$ or $s$ and $n \geq 2$ the Yoshida-lift $F = Y^{(n)}(\varphi_1, \varphi_2)$ is (if nonzero) an eigenfunction of the $p$-component of the Hecke-algebra of $GSp_n$ whose standard $L$-function is given for $n \geq 2$ by:*

$$\begin{aligned}
D_{F^{(n)}}^{(N)}(s) &= \prod_{p \nmid N} (1 - p^{-s})^{-1}(1 - \beta_p\widetilde{\beta}_p p^{-s})^{-1}(1 - \beta_p\widetilde{\beta}_p^{-1} p^{-s})^{-1}(1 - \beta_p^{-1}\widetilde{\beta}_p p^{-s})^{-1} \\
&\quad \cdot (1 - \beta_p^{-1}\widetilde{\beta}_p^{-1} p^{-s})^{-1} \prod_{j=1}^{n-2}(1 - p^{-s+j})^{-1}(1 - p^{-s-j})^{-1} \\
&= \zeta^{(N)}(s) L_{\varphi_1,\varphi_2}^{(N)}(s) \cdot \prod_{j=1}^{n-2} \zeta^{(N)}(s-j)\zeta^{(N)}(s+j).
\end{aligned}$$

*Here if $\varphi_1, \varphi_2$ are cuspidal and $\varphi_i$ corresponds to the elliptic cusp form*

$$f_i(z) = \sum_n a_i(n) \exp(2\pi i n z)$$

*under Eichler's correspondence we have*

$$L_{\varphi_1,\varphi_2}^{(N)}(s) = L_{f_1,f_2}^{(N)}\left(s + \frac{k_1 + k_2}{2} - 1\right),$$

*where $L_{f_1,f_2}^{(N)}(s)$ is the "good" part*

$$\zeta^{(N)}(2s + k_1 + k_2 - 2) \sum_{(n,N)=1} \frac{a_1(n)a_2(n)}{n^s}$$

*of the tensor product $L$-function of $f_1, f_2$.*

*Proof.* This is proved in the same way as Theorem 6.1 and Corollary 6.1 of [5]. □

**Theorem 4.2** *With the notations of Theorem 4.1 assume $\varphi_i$ to be essential for $i = 1, 2$. Then the function $\Lambda_N(s)$ of Section 2 is given by*

$$\Lambda_N(s) = \Lambda_N^{triv}(s) = \prod_{p|N} \prod_{j=1}^n (1 - p^{-s-2+j})^{-1}.$$

*If for some $p$ dividing $N$ only one of the $\varphi_i$ is $p$-essential (say $\varphi_1$), then there are precisely two local maximal orders $\tilde{R}_p$ and $\hat{R}_p$ strictly containing $R_p$. We let $\psi$ be the projection of $\varphi_2$ on the space of functions that are right invariant under the order $\tilde{R}$ with completion $\tilde{R}_p$ at $p$ and $R_l$ at the primes $l \neq p$, assume without loss of generality $\psi \neq 0$, and denote by $\alpha_p, \alpha_p^{-1}$ the $p$-Satake-parameters of $\psi$ and by $\epsilon_p$ the eigenvalue of $\varphi_2$ under $\tilde{w}_p$. Then the $p$-factor in the formula for $\Lambda_N(s)$ above has to be replaced by*

$$(1 + \epsilon_p \alpha_p p^{(-2s-1)/2})^{-1}(1 + \epsilon_p \alpha_p^{-1} p^{(-2s-1)/2})^{-1} \prod_{j=3}^n (1 - p^{-s-2+j})^{-1}$$

*Proof.* The proof is completely analogous to that of Corollary 7.1 of [5] (see also the note added after this corollary for the case that one of the $\varphi_i$ is not essential). There are two places where the modifications needed to cover the present case are not obvious. Firstly we have to generalize Evdokimov's [9] calculation of the action of the Hecke-operator

$$F \mapsto \sum_{\mathbf{A} \in \mathbf{Z}_{\text{sym}}^{(n,n)} \bmod p\mathbf{Z}_{\text{sym}}^{(n,n)}} F \Big| \begin{pmatrix} 1_n & \mathbf{A} \\ 0 & p1_n \end{pmatrix}.$$

This was used to prove the formula

$$p^{n(n+1)/2} \cdot s_p(S)^n \gamma_p^n \sum_{i=0}^{r/2} (-1)^i p^{i(i-1-n)} \sum_{j=1}^{h(i)} \frac{r(S, S_{m/2-i,j})}{e(S_{m/2-i,j})} \vartheta^{(n)}(S_{m/2-i,j})$$

for the action of the operator $K_p^{(n)} \begin{pmatrix} 0 & -p^{-1}1_n \\ p1_n & 0 \end{pmatrix} K_p^{(n)}$ (notation of [5]) on the theta series $\vartheta^{(n)}(S)$ (where we denote by $S_{ij}$ representatives of the classes of integral positive definite quadratic forms of rank $m$ and discriminant exactly divisible by $p^{2i}$ which are split over $\mathbf{Q}_p$, whose level is not divisible by $p^2$, and that are equivalent over all $\mathbf{Z}_\ell$ with $\ell \neq p$ to a fixed form $S$ (of square discriminant exactly divisible by $p^{m-r}$).

In the present situation we want to replace the theta series $\vartheta^{(n)}(S)$ by a theta series

$$\vartheta^{(n)}(L, q, P, Z) = \sum_{\mathbf{x} \in L^n} P(\mathbf{x}) \exp(2\pi i \text{tr}(q(\mathbf{x})Z))$$

with a $W_\rho$-valued harmonic polynomial $P$ as in Section 1, where the lattice $L$ has the Gram matrix $S$ from above with respect to the quadratic form $q$. A close examination of Evdokimov's proof shows then that the formula above remains true if one replaces the second sum by
$$\sum_{L_i} \vartheta^{(n)}(L_i, P, q, Z),$$
where $L_i$ runs over the lattices contained in $L$ whose Gram matrix with respect to $q$ is equivalent to one of the $S_{ij}$ from above and where the polynomial $P$ is (as in Section 1) considered as a function on the $n$-fold sum of the vector space that carries all the lattices involved.

We also need a generalization of Lemma 7.5 of [5] describing the commutation relation between Siegel's $\Phi$-operator and our Hecke operators for the bad places. We notice first that the formula given there should read:

$$F|T_n^{(n)}|\Phi = p^{2n-k} F|\Phi|T_{n-1}^{(n-1)},$$
$$F|T_i^{(n)}|\Phi = p^{n+i-k} F|\Phi|T_{i-1}^{n-1} + p^{2i} F|\Phi|T_i^{(n-1)} \ (i < n).$$

(In [5] the $k$ was replaced by 2). It is easy to see that this remains true if $F$ is a vector valued Siegel modular form transforming according to the representation of highest weight $(\mu_1, \ldots, \mu_{n-1}, k)$. In the case under consideration we have $k = 2$ (which is also true in the applications of this Lemma in [5]), and the argument proceeds as in loc.cit..

The rest of the proof goes through as in [5], replacing the ordinary theta series used there by the theta series with spherical harmonics used here in the same way as above (i. e. in formulas for the action of some Hecke operator on a theta series $\vartheta^{(n)}(L, q, P, Z)$ there appear theta series $\vartheta^{(n)}(L', q, P, Z)$ of lattices $L'$ on the same vector space $V$ and with the same polynomial $P$ on $V^n$).

## 5. Nonvanishing

In Section 1 we have reduced the question of the nonvanishing of Yoshida's lifting $Y^{(n)}(\varphi_1, \varphi_2, *)$ to the case where $*$ stands for $+$ or $s$ and where $n = 2$ and $\varphi_1, \varphi_2$ are non-proportional Hecke eigenforms. We can now put together the results of the previous section and formulate the answer to the nonvanishing problem in the following theorem:

**Theorem 5.1** *Let $\varphi_i \in \mathcal{A}(D_\mathbf{A}^\times, R_\mathbf{A}^\times, \tau_i)$ for $i = 1, 2$ be given and assume that for each $p \mid N$ at least one of $\varphi_1, \varphi_2$ is p-essential. Assume the $\varphi_i$ to be eigenfunctions of the of the involutions $\widetilde{w_p}$ for $p \mid N$ having the same eigenvalues and of the Hecke operators $\tilde{T}(p)$ for the $p \nmid N$. Assume moreover that the $\varphi_i$ are nonproportional and that they correspond to the normalized newforms $f, g$ of weight $k_i = 2 + 2\nu_i$ and levels $N_f, N_g$ with $\mathrm{lcm}(N_f, N_g) = N$ (and $\gcd(N_f, N_g)$ divisible by the ramified primes of* $D$) under Eichler's correspondence described in Section 1 (see also [7, Section 5] if they are cuspidal. Then for $*$ denoting $+$ or $s$ the Yoshida-lifting $Y^{(2)}(\varphi_1, \varphi_2, *)$ is zero if and only if $\nu_2 = 0$, one of $\varphi_1, \varphi_2$ is constant (without loss of generality $\varphi_2$), and $L(f, k_1/2) = 0$.

*Proof.* As shown in Section 1 the lifting $Y^{(2)}(\varphi_1, \varphi_2, *)$ is zero if and only if
$$Y^{(3)}(\varphi_1, \varphi_2, *)$$
is cuspidal. By Proposition 3.1 the latter implies that $\Lambda_N(s) D_F^{(N)}(s)$ is regular at $s = 1$. Theorems 4.1 and 4.2 give that this is equivalent to $L_{f,g}((k_1 + k_2)/2) = 0$ if both of $\varphi_1, \varphi_2$ are cuspidal. But by a result of Shahidi [24] this can not happen; an alternative proof for this in the setting of holomorphic modular forms could be obtained by generalizing Ogg's result [23] for the value at $s = 2$ of the tensor product $L$-function of two cusp forms of weight 2 to arbitrary weights. This proves the nonvanishing statement if $\varphi_1, \varphi_2$ are cuspidal. If $\varphi_1$ is cuspidal and $\varphi_2 = 1$ then $L_{\varphi_1, \varphi_2}^{(N)}(s) = L^{(N)}(f_1, s + \nu_1 + 1) L^{(N)}(f_1, s + \nu_1)$, and this is zero if and only if the central critical value $L(f_1, k_1/2)$ is zero. It remains to show that in this case one has indeed that $F = Y^{(2)}(\varphi_1, 1, *)$ is zero. For this we can again proceed in the same way as in [5], assuming $F \neq 0$. In order to see how to modify the argument let us sketch the argument given there: We used the following version of Siegel's theorem: The Eisenstein series $E_4^2(\mathcal{Z}, s)$ has a pole of order 1 in $s = 1/2$, and the residue $\mathcal{E}_4^2(\mathcal{Z})$ is a linear combination (with all coefficients $\alpha_i$ nonzero) of the genus theta series $\vartheta^{(4)}(\mathcal{G}_i)$ of degree 4 for the genera $\mathcal{G}_i$ of positive definite integral quadratic forms of rank 4, square discriminant and level dividing $N$. We then restricted $\mathcal{Z}$ to a block diagonal $\begin{pmatrix} Z_1 & 0 \\ 0 & Z_2 \end{pmatrix}$ and formed the double Petersson product $\langle\langle F, \mathcal{E}_4^2\rangle_{Z_1}, F\rangle_{Z_2}$. By Lemma 9.1 of [5] only the theta series of the genus of the order $R$ gave a contribution to this double product, which then became
$$c \sum_{i,j} \frac{|\langle F, \vartheta^{(2)}(I_{ij})\rangle|^2}{e_i e_j}$$
with $c \neq 0$. From this and Theorem 4.1 of [5] (or 2.29 of the present article) we then concluded that $\Lambda_N(s) D_F^{(N)}(s)$ has a pole in $s = 1$ unless $F$ is orthogonal to all the $\vartheta^{(2)}(I_{ij})$, which is excluded since $F$ by construction is a linear combination of these theta series. In our present situation we have to apply the differential operator $\Delta_2^\nu$ to all the $\vartheta^{(4)}(I_{ij})$ and obtain as in [1] after evaluating the double Petersson product
$$c \sum_{i,j,\kappa} \frac{|\langle F, \vartheta^{(2)}(I_{ij}, P_\kappa)\rangle|^2}{e_i e_j},$$
where $P_\kappa$ runs over a basis of the space of homogenous harmonic forms on $D_\infty^2$ transforming under the action of $GL_2(\mathbf{R})$ according to $\det^{2+\nu_1}$ (we notice that as in Lemma 9.1 of [5] $F$ is orthogonal to all the theta series with harmonic polynomials except those belonging to lattices in the genus of $R$). We conclude again that $\Lambda_N(s) D_F^{(N)}(s)$ has a pole in $s = 1$ unless $F$ is orthogonal to all the $\vartheta^{(2)}(I_{ij}, P_\kappa)$, which is excluded since $F$ by construction is a linear combination of these theta series with harmonic polynomials. By our formula for $D_F^{(N)}(s)$ in this situation from above we see that this implies $L(f_1, k_1/2) = 0$, which proves the rest of our assertion.

*Remark:* A different phrasing of Theorem 5.1 would say: Take any pair of cusp forms $f, g$ of weights $2 + 2\nu_i$ ($\nu_1 \geq \nu_2$) with trivial character for groups $\Gamma_0(N_{f,g})$ with square free $N_f, N_g$, $\gcd(N_f, N_g) \neq 1$. Assume $f, g$ to be newforms (in particular Hecke eigenforms) and that $f, g$ have the same eigenvalues under the Atkin-Lehner involutions $w_p$ for all $p \mid \gcd(N_f, N_g)$. Then we can construct a nonzero Siegel modular form of degree 2 transforming according to the representation of highest weight $(\nu_1 + \nu_2 + 2, \nu_1 - \nu_2 + 2)$ of $GL_2(\mathbf{C})$ under $\Gamma_0^{(2)}(N)$ with $N = \mathrm{lcm}(N_f, N_g)$, whose Satake parameters are related to those of $f, g$ as described in Section 4; this may be viewed as giving an explicit construction of an endoscopic lifting.

In order to construct this Siegel modular form we first replace $f, g$ by forms of level $N$ having the same Hecke eigenvalues as $f, g$ for $p \nmid N$ and being eigenfunctions of all the Atkin-Lehner involutions $w_p$ for $p \mid N$ with the same eigenvalues. We then fix a definite quaternion algebra over $\mathbf{Q}$ whose ramified primes divide $\gcd(N_f, N_g)$ and use the generalized correspondence of Eichler for an Eichler order of level $N$ in $D$ from Section 1 and [7, Section 5] to find $\varphi_1, \varphi_2$ corresponding to $f, g$ as required in Theorem 5.1 and apply this theorem to the pair $\varphi_1, \varphi_2$. It should be noted that when $\gcd(N_f, N_g)$ is not prime this allows several possible constructions whose results are pairwise orthogonal by the generalization of Lemma 9.1 of [5] noted in the proof of theorem 5.1. Such a possibility has been discussed under the name of twin forms in [15].

In particular for $\nu_1 = 1$, $\nu_2 = 0$ we obtain scalar valued Siegel modular forms of weight 3 and degree 2, which in view of their connection to differential forms on the quotient of $\mathbf{H}_2$ by the group $\Gamma_0^{(2)}(\mathrm{lcm}(N_1, N_2))$ are of interest in geometry.

## 6. An example

To illustrate our results we want to discuss an example. The calculations for this example were done using the algorithms from [26] and the PARI number theory package. We put $N = 17$. The quaternion algebra over $\mathbf{Q}$ that is ramified at $\infty$ and at $N$ only has type number and class number $t = h = 2$. A basis is given by the $f_i^{(0)}$ for $i = 0, \ldots, 3$ with $f_0^{(0)} = 1$ and $(f_1^{(0)})^2 = -17$, $(f_2^{(0)})^2 = -3$, $f_1^{(0)} f_2^{(0)} = -f_2^{(0)} f_1^{(0)} = f_3^{(0)}$. The $\mathbf{Z}$-lattice generated by these vectors has index 12 in a maximal order. It is transformed into a maximal order $R_0$ containing it by the basis transformation with

matrix
$$\begin{pmatrix} 1 & \frac{1}{2} & 0 & 0 \\ 0 & 0 & \frac{1}{2} & 0 \\ 0 & \frac{1}{2} & 0 & \frac{1}{3} \\ 0 & 0 & \frac{1}{2} & \frac{-1}{3} \end{pmatrix}$$

and into another maximal order $R_1$ by the basis transformation with matrix
$$\begin{pmatrix} 1 & \frac{1}{2} & \frac{1}{2} & 0 \\ 0 & \frac{1}{4} & \frac{1}{4} & \frac{1}{4} \\ 0 & \frac{-1}{3} & \frac{3}{4} & \frac{-1}{2} \\ 0 & \frac{1}{12} & \frac{-1}{8} & \frac{-1}{4} \end{pmatrix}.$$

The resulting lattice has Gram matrix (with respect to $\mathrm{tr}(x\bar{y})$ and to the new basis $\{f_i\}$)
$$A_1 = \begin{pmatrix} 2 & 1 & 1 & 0 \\ 1 & 4 & -1 & 1 \\ 1 & -1 & 6 & 2 \\ 0 & 1 & 2 & 10 \end{pmatrix}.$$

The product in $R_1$ of $\sum_{i=0}^{3} \alpha_i f_i$ and $\sum_{i=0}^{3} \beta_i f_i$ is $\sum_{i=0}^{3} \gamma_i f_i$ with

$$\begin{aligned}
\gamma_0 &= \alpha_0\beta_0 + ((-2\alpha_1 + (\alpha_2+\alpha_3))\beta_1 + (-2\beta_3\alpha_1 + ((-3\alpha_2 - 3\alpha_3)\beta_2 + (\beta_3\alpha_2 - 5\alpha_3\beta_3)))) \\
\gamma_1 &= \alpha_1\beta_0 + (\beta_1\alpha_0 + ((\alpha_1 - \alpha_3)\beta_1 + ((\beta_2+\beta_3)\alpha_1 + (3\alpha_3\beta_2 - 3\beta_3\alpha_2)))) \\
\gamma_2 &= \alpha_2\beta_0 + (\beta_2\alpha_0 + (-2\alpha_3\beta_1 + ((\beta_2+2\beta_3)\alpha_1 + ((\alpha_2+\alpha_3)\beta_2 - \beta_3\alpha_2)))) \\
\gamma_3 &= \alpha_3\beta_0 + (\beta_3\alpha_0 + ((\alpha_2+\alpha_3)\beta_1 + (-\beta_2\alpha_1 + \beta_3\alpha_2))).
\end{aligned}$$

The other type of maximal order is $R_2$, which is obtained from $R_1$ by the basis transformation with matrix
$$\begin{pmatrix} \frac{3}{4} & 1 & \frac{1}{2} & 0 \\ 0 & 0 & -1 & 1 \\ \frac{-1}{2} & 0 & -1 & 0 \\ \frac{1}{4} & 0 & \frac{-1}{2} & -1 \end{pmatrix},$$

it has Gram matrix
$$\begin{pmatrix} 2 & 1 & 0 & 0 \\ 1 & 2 & -1 & 1 \\ 0 & -1 & 12 & 5 \\ 0 & 1 & 5 & 12 \end{pmatrix}.$$

An ideal $I_{12}$ with left order $R_2$ and right order $R_1$ is obtained from $R_1$ by the basis transformation with matrix
$$\begin{pmatrix} 0 & 1 & -1 & \frac{1}{2} \\ \frac{1}{2} & -1 & \frac{-1}{2} & -1 \\ 0 & 0 & 0 & -1 \\ \frac{1}{2} & 0 & \frac{1}{2} & \frac{1}{2} \end{pmatrix},$$

it has Gram matrix
$$\begin{pmatrix} 4 & -2 & 1 & 0 \\ -2 & 4 & 0 & 1 \\ 1 & 0 & 6 & 3 \\ 0 & 1 & 3 & 6 \end{pmatrix}$$
with respect to the basis $\{f'_i\}$ obtained. The order $R_2$ has 6 units, the order $R_1$ has no nontrivial units. We obtain therefore as the unique (up to scalar multiples) nontrivial function in $\mathcal{A}(D_{\mathbf{A}}^\times, R_{\mathbf{A}}^\times, 1)$ orthogonal to the constants the function $\varphi_2$ with $\varphi_2(y_1) = -1$, $\varphi_2(y_2) = 3$ (corresponding to the unique normalized cusp form of weight 2 and level 17). The space of cusp forms of weight 4 and level 17 has dimension 4. It is easily checked by looking at the linear forms on $D \otimes \mathbf{R}$ and using Eichler's correspondence described in Section 1 that only one of these has $w_{17}$-eigenvalue $-1$. The corresponding function $\varphi_1 \in \mathcal{A}(D_{\mathbf{A}}^\times, R_{\mathbf{A}}^\times, \tau_1)$ is given by $\varphi_2(y_2) = 0$ and $\varphi_2(y_1)(\sum_{i=0}^3 \alpha_i f_i) = \alpha_3$. Using this we find that the Siegel cusp form $Y^{(2)}(\varphi_1, \varphi_2)$ of weight 3 is given as the sum of the theta series of degree 2 of $R_1$ with polynomial

$$P_1(\sum_{i=0}^3 \alpha_i f_i, \sum_{i=0}^3 \beta_i f_i) = 2\alpha_3\beta_0 + (-2\beta_3\alpha_0 + ((2\alpha_2 + 2\alpha_3)\beta_1 + (-2\beta_2 - 2\beta_3)\alpha_1))$$

and the theta series of degree 2 of $I_{12}$ with polynomial

$$P_{12}((\sum_{i=0}^3 \alpha'_i f'_i, \sum_{i=0}^3 \beta'_i f'_i)) = -2\alpha_2\beta_0 + (2\beta_2\alpha_0 + (-2\alpha_3\beta_1 + 2\beta_3\alpha_1)).$$

The computation gives the following beginning terms of the Fourier expansion of $Y^{(2)}(\varphi_1, \varphi_2)$ (where the variable in $\mathbf{H}_2$ is
$$Z = \begin{pmatrix} \tau & z \\ z & \tau' \end{pmatrix}$$
and where we put $X_1 = \exp(2\pi i \tau)$, $X_2 = \exp(2\pi i z)$, $X_3 = \exp(2\pi i \tau'))$ so that the coefficient at $X_1^a X_2^b X_3^c$ is the Fourier coefficient at $\begin{pmatrix} a & b/2 \\ b/2 & c \end{pmatrix}$):

$$(-32X_3^6 X_2^2 - 64X_3^6 X_2)X_1^5$$
$$+(-32X_3^5 X_2^3 - 96X_3^6 X_2^2 + (32X_3^6 - 64X_3^5)X_2)X_1^4$$
$$+((-32X_3^6 + 32X_3^5 + 32X_3^4)X_2^2 - 32X_3^6 X_2)X_1^3$$
$$+(-32X_3^5 - 32X_3^4 + 32X_3^3)X_2 X_1^2.$$

(We computed the coefficients of $X_1^a X_2^b X_3^c$ for all reduced binary forms $[a, b, c] = ax^2 + bxy + cy^2$ with $c \le 25$) To be sure that this is indeed a Hecke eigenform we computed the action of $T(2)$, $T(3)$ and $T(5)$ following the formulas given in [27]; it turns out that the Hecke eigenvalues are $-5$, $-8$ and $-4$ respectively. Notice that in the formula at the bottom of p. 386 of [27] the term $a([m/p, r, np])$ should be replaced by $a([m/p, -r, np])$ (which matters since we are dealing with odd weight).

Siegfried Böcherer, Fakultät für Mathematik und Informatik, Universität Mannheim, Seminargebäude A5, D-68131, Mannheim, Germany

Rainer Schulze-Pillot, Fachbereich Mathematik, Universität des Saarlandes, Postfach 151150, D-66041 Saarbrücken, Germany